\documentstyle[amssymb,amsfonts,12pt]{article}

\newtheorem{corollary}{Corollary}

\newtheorem{lemma}{Lemma}
\newtheorem{definition}{Definition}
\newtheorem{theorem}{Theorem}

\newtheorem{claim}{Claim}
\newtheorem{subclaim}{Subclaim}

\newtheorem{fact}{Fact}

\newcommand{\bds}{\begin{description}}
\newcommand{\eds}{\end{description}}
\newcommand{\be}{\begin{enumerate}}
\newcommand{\ee}{\end{enumerate}}
\newcommand{\bq}{\begin{quote}}
\newcommand{\eq}{\end{quote}}

\newcommand{\nothing}[1]{}

\def\cal{\mathcal}
\newcommand{\MM}{{\sf MM}}
\newcommand{\PFA}{{\sf PFA}}
\newcommand{\BPFA}{{\sf BPFA}}

\newcommand{\ZFC}{{\sf ZFC}}

\newcommand{\SCH}{\sf SCH}
\newcommand{\GCH}{\sf GCH}

\newcommand{\w}{\omega}

\newcommand{\MRP}{{\sf MRP}}

\newcommand{\cf}{{\mbox{cof}}}

\renewcommand{\colon}{\,\mbox{\tt:}\,}

\begin{document}

\title{The Proper Forcing Axiom and the Singular Cardinal
Hypothesis}\footnotetext{\hspace{2mm} 2000 Mathematical Subject
Classification: 03E05, 03E10, 03E65, 03E75. Key words: $\MRP$,
$\PFA$, $\SCH$.}

\author
{Matteo Viale \thanks{\hspace{2mm} Dipartimento di Matematica,
Universit\'a di Torino and Equipe de Logique Math\'ematique,
Universit\'e Paris 7, email: {\texttt viale@dm.unito.it}}}

\date{}

\maketitle

\abstract{We show that the Proper Forcing Axiom implies the Singular Cardinal Hypothesis.
The proof uses the reflection principle $\MRP$ introduced by Moore in \cite{mooMRP}. }

\section*{Introduction}

In one of the first applications of the forcing techniques, Easton
\cite{eas63} showed that the exponential function $\kappa \mapsto
2^\kappa$ on regular cardinals can be arbitrary modulo some mild
restrictions.  The situation for singular cardinals is much more
subtle.  For instance, Silver \cite{sil75} showed that the
Singular Cardinal Hypothesis $\SCH$ cannot first fail at a
singular cardinal of uncountable cofinality. Recall that $\SCH$
states that $2^{\kappa}=\kappa^+$, for all singular strong limit
cardinals $\kappa$.  However, it is known that $\SCH$ can first
fail even at $\aleph_{\omega}$.

The role of large cardinals in this context is twofold. On one
hand they are necessary for the construction of models of the
negation of $\SCH$ since any such model has an inner model with measurable cardinals (see \cite{gitSCH} for a survey of Prikry type
forcings and applications to $\SCH$). On the other hand it is a
theorem of Solovay \cite{solSCH} that $\SCH$ holds above a
strongly compact cardinal. Forcing axioms imply reflection
principles similar to the one used in Solovay's proof, thus it was
reasonable to expect that they would also settle $\SCH$. Indeed,
in \cite{magforshe}, Foreman, Magidor and Shelah showed that the
strongest forcing axiom, Martin's Maximum $\MM$, implies $\SCH$.
This was later improved by Veli\v{c}kovi\'c \cite{velSCH} who also
showed that $\SCH$ follows from $\PFA ^+$. In fact, what is shown
in \cite{velSCH} is that if $\theta
>\aleph_1$ is regular and stationary subsets of
$[\theta]^{\omega}$ reflect to an internally closed and unbounded
set, then $\theta ^{\omega}=\theta$. This, combined with Silver's
theorem, implies $\SCH$. At this point, it was left open whether
$\SCH$ is a consequence of $\PFA$.

Very little progress was made on this problem for over fifteen
years. Then, in 2003, Moore  \cite{mooMRP} introduced a new
reflection principle, the {\it Mapping Reflection Principle}
$\MRP$ and deduced it from $\PFA$. He showed that $\MRP$ implies
the continuum is equal to $\aleph_2$ and the failure of
$\Box({\kappa})$, for all $\kappa >\aleph_1$. $\MRP$ has many
features in common with the reflection principles which follow
from $\MM$, so it should be expected that $\MRP$ could affect the
behaviour of the exponential function also on higher cardinals. In
fact, Moore showed in \cite{mooSCH} that if $\MRP$ holds and
$\kappa >\omega_1$ is a regular cardinal with a nonreflecting
stationary set consisting of points of countable cofinality, then
$\kappa^{\omega_1}=\kappa$. This, combined with the above result
of Veli\v{c}kovi\'c, strongly suggests that $\PFA$ implies $\SCH$.
In this paper we confirm this conjecture.

The paper is organized as follows.
In section \S \ref{mainsec} we prove that if $\MRP$ holds then
$\kappa^\omega=\kappa$ for all regular cardinals $\kappa
>\omega_1$. This, together with Silver's theorem, implies that
$\SCH$ follows from $\PFA$. In the final section we show that the result can be refined
and applied to study another interesting problem in the area of forcing axioms, i.e. to
investigate what kind of forcing notions can preserve this type of axioms.

Our notation is standard and follows \cite{jecST} and
\cite{mooMRP}. For a regular cardinal $\theta$, we use $H(\theta)$
to denote the structure $(H(\theta),\in,<)$  whose domain is the
collection of sets whose transitive closure is of size less than
$\theta$ and where $<$ is a predicate for a fixed well ordering of
$H(\theta)$. If $X$ is an uncountable set, $C\subseteq[X]^\omega$
is closed and unbounded (club) if there is
$f:[X]^{<\omega}\rightarrow X$ such that $C$ is the set of all
$Y\in [X]^{\omega}$ such that $f[Y]^{<\omega}\subseteq Y$.
$S\subseteq[X]^\omega$ is  stationary if it intersects all club
subsets of $[X]^{\omega}$. The $f$-closure of $X$ is the smallest
$Y$ containing $X$ such that $f[Y]^{<\omega}\subseteq Y$. If $X$
is a set of ordinals then $\overline{X}$ denotes the topological
closure of $X$ in the order topology.


\section{The main result} \label{mainsec}

The purpose of this section is to show that $\MRP$ implies that
$\lambda^\omega=\lambda$, for every cardinal $\lambda > \omega_1$
with $\cf(\lambda)>\omega$. We start by recalling the relevant
definitions from \cite{mooMRP}.

\begin{definition}
Let $\theta$ be a regular cardinal, let $X$ be uncountable, and
let $M\prec H(\theta)$ be countable such that $[X]^{\omega}\in M$.
A subset $\Sigma$ of $[X]^{\omega}$ is $M$-stationary if
for all $E\in M$ such that $E\subseteq [X]^{\omega}$ is club,
$\Sigma \cap E\cap M \neq \emptyset$.
\end{definition}

Recall that the Ellentuck topology on $[X]^{\omega}$ is obtained
by declaring a set open if it is the union of sets of the form
$$
[x,N]=\{ Y\in [X]^{\omega} \colon x\subseteq Y\subseteq N \}
$$
where $N\in[X]^\w$ and $x\subseteq N$ is finite.

\begin{definition}
$\Sigma$ is an open stationary set mapping if there is an
uncountable set $X$ and a regular cardinal $\theta$
such that $[X]^{\omega}\in H(\theta)$, the domain of $\Sigma$ is a
club in $[H(\theta)]^{\omega}$ of countable elementary submodels
$M$ such that $X\in M$ and for all $M$, $\Sigma(M)\subseteq
[X]^{\omega}$ is open in the Ellentuck topology on $[X]^\omega$ and $M$-stationary.
\end{definition}

\noindent The Mapping Reflection Principle ($\MRP$) asserts that:
\begin{quote}
If $\Sigma$ is an open stationary set mapping, there is a continuous $\in$-chain $\vec N=
(N_\xi\colon\xi<\w_1)$ of elements in the domain of $\Sigma$ such
that for all limit ordinals $\xi<\w_1$ there is  $\nu<\xi$ such
that $N_\eta\cap X\in\Sigma(N_\xi)$ for all $\eta$ such
that $\nu <\eta <\xi$.
\end{quote}

\noindent If $(N_{\xi}\colon \xi <\omega_1)$ satisfies the
conclusion of $\MRP$ for $\Sigma$ then it is said to be a
reflecting sequence for $\Sigma$. It is shown in \cite{mooMRP}
that $\MRP$ is a consequence of $\PFA$. We are now ready to prove
the following theorem.

\begin{theorem} \label{mainth}
Assume $\MRP$. Then $\lambda^{\aleph_0}=\lambda$, for every
$\lambda \geq \omega_2$ of uncountable cofinality.
\end{theorem}

\noindent {\bf Proof:} We will prove the theorem by induction. The
base case $\lambda =\aleph_2$ is handled by Moore's result
\cite{mooMRP} that $\MRP$ implies $2^{\aleph_0}=\aleph_2$. If
$\lambda =\kappa ^+$ with $\cf(\kappa)> \omega$ then
$\lambda^{\aleph_0}=\lambda \cdot \kappa ^{\aleph_0}$, so the
result holds by the inductive hypothesis. If $\lambda$ is a limit
cardinal and $\cf(\lambda)>\omega$ then $\lambda^{\aleph_0}= \sup
\{ \mu^{\aleph_0} : \mu <\lambda\}$, so the result also follows by
the inductive hypothesis. Thus, the only interesting case is when
$\lambda=\kappa^+$, with $\kappa$ singular of countable
cofinality. In this case we will show, using $\MRP$, that
$\kappa^{\aleph_0}=\kappa^+$.

Now, let $\kappa$ be singular of countable cofinality and assume
the theorem holds below $\kappa$. Fix a sequence
$(C_\delta\,:\,\delta\in\kappa^+)$ such that $C_{\delta}$ is a
club in $\delta$ of minimal order type. In fact, we will be
interested only in ordinals $\delta$ of cofinality $\leq\omega_1$. For every pair of ordinals $\delta, \beta <\kappa^+$,
we fix a decomposition $\delta=\bigcup_n K(n,\delta,\beta)$ such
that:

\bds

\item[\it (i)] $|K(n,\delta,\beta)|<\kappa$, for all
$n$

\item[\it (ii)]   $K(n,\delta,\beta)\subseteq K(m,\delta,\beta)$,
for $n<m$

\item[\it (iii)]  if $\eta<\beta$ is of cofinality $\omega_1$ then there is $n$ such that $C_\eta\cap\delta\subseteq K(n,\delta,\beta)$

\item[\it (iv)]  $K(n,\delta,\beta)$ is a closed subset of $\delta$, for all $n$.

\eds

\noindent This is easily achieved, for example, as follows. First
of all, fix
 an increasing sequence $(\kappa_n:n\in\omega)$ of regular cardinals
converging to $\kappa$. For all $\eta<\kappa^+$ let
$\phi_\eta:\kappa\rightarrow\eta$ be a surjection. Now set:
$$
K(n,\delta,\beta)=\delta\cap\overline{\phi_\delta[\kappa_n]\cup\bigcup\{C_\eta\cap\delta\,:\eta\in\phi_\beta[\kappa_n]\,\&\,\cf\,\eta=\omega_1\}}.
$$

\noindent Fix also a partition $\{A_s\,:\,s\in\kappa^{<\omega}\}$
 of $\{ \delta <\kappa^+: \cf(\delta)=\omega \}$ into disjoint stationary
sets. Let $D(n,\delta,\beta)$ be the set of all
$g\in\kappa^\omega$ such that there are infinitely many $j$ such
that $K(n,\delta,\beta)\cap A_{g\upharpoonright j} \neq
\emptyset$. Using the fact that $K(n,\delta,\beta)$ is of size
$<\kappa$ and the inductive hypothesis we immediately have the
following.

\begin{fact} \label{boundfa}
$D(n,\delta,\beta)$ is of size smaller than $\kappa$, for all
$n,\delta$ and $\beta$. \hfill$\square$\medskip
\end{fact}

\noindent We will be done once we show the following.

\begin{lemma} Assume $\MRP$. Then
$\bigcup \{ D(n,\delta,\beta) : n<\omega \mbox{ and } \delta,
\beta <\kappa^+ \}=\kappa^{\omega}$.
\end{lemma}

\noindent {\bf Proof:}  Fix $g\in\kappa^\omega$. We have to find
some $(n,\delta,\beta)$ such that $g\in D(n,\delta,\beta)$. We are
going to define an open stationary set mapping $\Sigma_g$ and
apply $\MRP$. We first fix some notation. Given a countable set
$X$, we let $\delta_X=\sup(X\cap\kappa^+)$ and
$\alpha_X=\sup(X\cap\omega_1)$. If $\alpha<\gamma<\omega_1$, let
the {\em height} of $\alpha$ in $\gamma$ be defined by
$ht_\gamma(\alpha)=|C_\gamma\cap\alpha|$. Fix a sufficiently large
regular cardinal $\theta$. Suppose $M$ is a countable elementary
submodel of $H(\theta)$ containing all the relevant information.
Fix $\beta_M< \kappa^+$ large enough such that for every
$\gamma<\kappa^+$ of cofinality $\omega_1$ there is $\eta<\beta_M$
of cofinality $\omega_1$ such that $C_\gamma\cap M=C_\eta\cap M$.
If $\gamma\in M\cap\kappa^+$ let  $n_\gamma$ be the smallest
integer $l$ such that $\gamma\in K(l,\delta_M,\beta_M)$. Now, let
$\Sigma_g(M)$ be the set of all $X\in[M\cap\kappa^+]^\omega$ such
that $\alpha_X<\alpha_M$, $\delta_X<\delta_M$ and such that letting $m=ht_{\alpha_M}(\alpha_X)$,
we have that $A_{g\upharpoonright m}\cap
K(n_{\delta_X},\delta_M,\beta_M)\neq\emptyset$. We will show that
$\Sigma_g(M)$ is open and $M$-stationary, for all $M$.

\begin{claim} $\Sigma_g(M)$ is open.
\end{claim}

\noindent {\bf Proof:} Suppose $X\in\Sigma_g(M)$. First find
$\alpha\in X\cap\omega_1$ such that
$ht_{\alpha_M}(\alpha)=ht_{\alpha_M}(\alpha_X)=m$.  Let
$n=n_{\delta_X}$. By the definition of $n_{\delta_X}$ we have that
$\delta_X\not\in K(n-1,\delta_M,\beta_M)$. Since
$K(n-1,\delta_M,\beta_M)$ is a closed subset of $\delta_M$, there
is a $\gamma$ below $\delta_X$ such that $(\gamma,\delta_X]\cap
K(n-1,\delta_M,\beta_M)=\emptyset$. Pick $\delta\in X$ larger than
$\gamma$. Then we have the following.

\begin{subclaim}
$[\{\alpha,\delta\},X]\subseteq\Sigma_g(M)$.
\end{subclaim}

\noindent {\bf Proof:} If $Y\in[\{\alpha,\delta\},X]$, clearly, we
have that $ht_{\alpha_M}(\alpha_Y)=m$. Since $\delta \in Y$ and
$Y\subseteq X$ we have that $\delta \leq \delta_Y\leq \delta_X$.
By the above remarks we can conclude that $n=n_{\delta_X}\leq
n_{\delta_Y}$. Now, since
$X\in\Sigma_g(M)$,
$K(n_{\delta_X},\delta_M,\beta_M)\cap A_{g\upharpoonright m}
\neq\emptyset$.
By construction,
$K(n_{\delta_X},\delta_M,\beta_M)\subseteq K(n_{\delta_Y},\delta_M,\beta_M)$ so
$K(n_{\delta_Y},\delta_M,\beta_M)\cap A_{g\upharpoonright m}
\neq\emptyset$.
Since $ht_{\alpha_M}(\alpha_Y)=m$, $Y\in
\Sigma_g(M)$. \hfill$\square$\medskip

\begin{claim} $\Sigma_g(M)$ is $M$-stationary.
\end{claim}
\noindent {\bf Proof:} Given
$f:[\kappa^+]^{<\omega}\rightarrow\kappa^+$ belonging to $M$, we
must find $X\in M\cap\Sigma_g(M)$ which is closed under $f$.
First, find $N\in M$, a countable elementary submodel of
$H(\kappa^{++})$ containing all the relevant objects for the
argument below. Let $ht_{\alpha_M}(\alpha_{N})=m$ and  find
$\alpha\in N$ with the same height in $\alpha_M$. Now let $C$ be
the set of $\delta <\kappa^+$ such that
$f[\delta]^{<\omega}\subseteq \delta$. Then $C$ is a club subset
of $\kappa^+$ and $C\in N$. Since $A_{g\restriction m}$ is
stationary in $\kappa^+$ and, by our assumption it belongs to $N$,
we can find $\delta\in C\cap A_{g\upharpoonright m}\cap N$. Then
$\delta\in K({n_\delta},\delta_M,\beta_M)\cap A_{g\restriction m}$.
Finally, let $Z\in N$ be a countable set cofinal in
$\delta$ and let $X$ be the $f$-closure of $\{\alpha\}\cup Z$.
Then $\delta_X=\delta$ and
$$
m=ht_{\alpha_M}(\alpha)\leq ht_{\alpha_M}(\alpha_X)\leq
ht_{\alpha_M}(\alpha_{N})= m.
$$
Since $K(n_{\delta_X},\delta_M,\beta_M)\cap
A_{g\upharpoonright m}$ is nonempty,
$X\in\Sigma_g(M)\cap M$ and $X$ is closed under $f$.
\hfill$\square$\medskip

Let $(M_\eta\,:\, \eta<\omega_1)$ be  a reflecting sequence for
$\Sigma_g$ provided by $\MRP$. Let $N=\bigcup_\eta M_\eta$ and
$\delta=\sup ( N\cap\kappa^+)$. Let $\delta_{\eta}=\sup(M_{\eta}\cap
\kappa^+)$, for every $\eta <\omega_1$.  We find a club
$E\subseteq \omega_1$ such that $\{\delta_\eta : \eta \in E
\}\subseteq C_\delta$ and $M_{\eta}\cap \omega_1=\eta$, for all
$\eta \in E$. Let $\alpha$ be a limit point of $E$.  For the rest
of this proof let $M=M_\alpha$. Now $C_\delta\cap M=C_\gamma\cap
M$ for some $\gamma<\beta_M$, by the choice of $\beta_M$. By condition {\it
(iii)} on $K(i,\delta_M,\beta_M)$ there is an
$n$ such that $C_\gamma\cap\delta_M$ is a subset of
$K(n,\delta_M,\beta_M)$. Since $C_\gamma\cap M\subseteq
C_\gamma\cap\delta_M$, we can conclude that $C_\delta\cap M$ is a
subset of $K(n,\delta_M,\beta_M)$.

Let $\nu<\alpha$ be such that $M_\eta\in\Sigma_g(M)$, for all
$\eta$ such that $\nu < \eta < \alpha$. For any such $\eta\in E$,
$M_\eta\in M$, so $\delta_\eta\in C_\delta\cap M\subseteq
K(n,\delta_M,\beta_M)$. If $\eta \in E$ and
$ht_{\alpha_M}(\eta)=j$, then $A_{g\upharpoonright j}\cap
K(n_{\delta_\eta},\delta_M,\beta_M)\neq\emptyset$ and, since
$n_{\delta_\eta}\leq n$, we have that
$K(n_{\delta_\eta},\delta_M,\beta_M)\subseteq
K(n,\delta_M,\beta_M)$. Now, for any $i$ we can find an $\eta\in
E$ such that $\nu <\eta <\alpha$ and $ht_{\alpha_M}(\eta)\geq i$,
so there are infinitely many $j$ such that $A_{g\upharpoonright
j}\cap K(n,\delta_M,\beta_M)\neq\emptyset$, i.e. $g \in
D(n,\delta_M,\beta_M)$, as desired. \hfill$\blacksquare$\medskip

\noindent Now we have the following immediate corollary.

\begin{corollary}$\PFA$ implies $\SCH$.
\end{corollary}

\noindent {\bf Proof:} This follows by induction. By Silver's
theorem \cite{sil75} the first cardinal violating $\SCH$ cannot be
singular strong limit of uncountable cofinality. On the other
hand, if $\kappa$ is a singular strong limit cardinal of countable
cofinality then, by Theorem 1,
$2^{\kappa}=\kappa^{\aleph_0}=\kappa^+$.
\hfill$\blacksquare$\medskip

\bigskip

\section{Final remarks and side results}

The techniques presented in the previous sections can be applied
to investigate another interesting problem in the area of forcing
axioms. Since forcing axioms have been able to settle many of the
classical problems of set theory, we can expect that the models of
a forcing axiom are in some  sense canonical. There are many
ways in which one can give a precise formulation to this concept.
For example, one can study what kind of forcings can preserve
$\PFA$, or else if a model $V$ of a forcing axiom can have an
interesting inner model $M$ of the same forcing axiom. There are
many results in this area, some of them very recent. For instance,
 K\"onig and Yoshinobu {\cite[Theorem 6.1]{KY}} showed that
$\PFA$ is preserved by $\w_2$-closed forcing. The same holds for
$\BPFA$. In fact, $\BPFA$ is preserved by any proper forcing that
does not add subsets of $\w_1$. In the other direction, in
\cite{velSCH} Veli\v{c}kovi\'c showed that if $\MM$ holds and $M$
is an inner model such that $\w_2^M=\w_2$, then ${\cal
P}(\w_1)\subseteq M$ and  in a very recent paper Caicedo and
Veli\v{c}kovi\'c \cite{velcaiBPFA} showed that if
$M\subseteq V$ are models of $\BPFA$ and $\omega_2^M=\omega_2$
then ${\cal P}(\w_1)\subseteq M$. Their argument also shows that if
$M\subseteq V$ are models of $\MRP$ and $\omega_2^M=\omega_2$, then
${\cal P}(\w)\subseteq M$.
We can use the result of the
previous section combined with this last result to show that
$\PFA$ is destroyed by many of the cardinal preserving notions of
forcing which add new $\omega$-sequences. A result of this sort
has been obtained by Moore in \cite{mooSCH}.

\begin{theorem}
Let $V$ and $W$ be two models of set theory with the same
cardinals with $V\subseteq W$. Assume $V$ and $W$ are both models
of $\MRP$  and that, moreover, for every cardinal $\kappa$,
stationary subsets in $V$ of $\{\alpha\in\kappa^+ :
 \, \cf\,\alpha=\omega\}$ are also stationary in $W$.
 Then $V$ and $W$ have the same $\omega$-sequences of ordinals.
\end{theorem}

\noindent {\bf Proof:} Assume otherwise. We proceed by induction
on the least cardinal $\kappa$ such that there in an
$\omega$-sequence of elements of $\kappa$ which is in $W$, but not
in $V$.  The base case $\omega$ is handled by the above result of
Caicedo and Veli\v{c}kovi\'c. We now run into two cases: either
the least such $\kappa$ has countable cofinality in $V$ or it
doesn't. The more involved case appears when
$\cf^V\,\kappa>\omega$. We present in some detail how to prove the
induction step in this situation. With minor modifications the
reader can supply the proof for the case that $\kappa$ is of
countable cofinality. The idea is to redo the proof of the
previous section using a $g\in\kappa^\omega\setminus V$. However
some extra care has to be paid in the definition of the sets
$K(n,\delta,\beta)$. Let $\{A_s:s\in\kappa^{<\omega}\}\in V$ be a
partition of the points of countable $V$-cofinality  of $\kappa^+$
into disjoint stationary sets. By the assumptions each $A_s$ is
still stationary in $W$. Fix
$(E_\delta:\delta<\kappa^+\,\&\,\cf^V\,\delta\geq\omega_1)\in V$
such that for all $\delta$ in its domain, $E_\delta$ is a club in
$\delta$ of minimal $V$-order-type.  So for each $\delta$,
$E_\delta$ has order type at most $\kappa$. Define in $V$, for all
$\alpha<\kappa$ and $\delta,\beta<\kappa^+$, sets
$K(\alpha,\delta,\beta)$ such that $\delta =\bigcup_\alpha
K(\alpha,\delta,\beta)$ and:

\bds

\item[\it (i)] $|K(\alpha,\delta,\beta)|<\kappa$

\item[\it (ii)]   $K(\alpha,\delta,\beta)\subseteq K(\gamma,\delta,\beta)$
for $\alpha<\gamma$

\item[\it (iii)]  if $\eta<\beta$ and $\cf^V\,\eta\geq\omega_1$ then there is $\alpha$ such that $E_\eta\cap\delta\subseteq K(\alpha,\delta,\beta)$

\item[\it (iv)]  $K(\alpha,\delta,\beta)$ is a closed subset of $\delta$.

\eds

\noindent This is easily achieved, for example, as follows. For all $\eta\in[\kappa,\kappa^+)$ let
$\phi_\eta:\kappa\rightarrow\eta$ be a bijection. Now set:
$$
K(\alpha,\delta,\beta)=\delta\cap\overline{\phi_\delta[\alpha]\cup\bigcup
\{E_\eta\cap\delta\,:\eta\in\phi_\beta[\alpha]\: \&\:
|E_\eta\cap\delta|\leq|\alpha|\}}.
$$

\noindent Define $D(\alpha,\delta,\beta)$ to be the set of all
$g\in\kappa^\omega$ such that there are infinitely many $j$ such
that $A_{g\upharpoonright j}\cap
K(\alpha,\delta,\beta)\neq\emptyset$ and use the inductive
hypothesis to get that
$D(\alpha,\delta,\beta)^V=D(\alpha,\delta,\beta)^W$. Since
$\kappa$ is the least cardinal with a new $\omega$-sequence, it
follows that $\cf^{W}(\kappa)=\omega$. Let
$g=(\alpha_n:n\in\omega)\in W$ be cofinal in $\kappa$. From now on
work in $W$. Let $K(n,\delta,\beta):=K(\alpha_n,\delta,\beta)$.
Now as in the previous section
use the parameters
$$
\{K(n,\delta,\beta):n<\omega\,\&\,\delta,\beta<\kappa^+\},\,
\{A_s:s\in\kappa^{<\omega}\},\,g
$$
\noindent to define $\Sigma_g$ and show that it is an open
stationary set mapping. We also refer to the previous section for
the notation.

Now, apply $\MRP$ to $\Sigma_g$ and using exactly the same argument as in the previous section show that
$g\in D(\alpha_n,,\delta_M,\beta_M)$ for some $M$ in a reflecting sequence for $\Sigma_g$. Some extra care has to be paid since if $(M_\eta\,:\,\eta<\omega_1)$ is  a reflecting sequence provided by $\MRP$ and
$\delta=\sup_{\eta<\omega_1} \delta_{M_\eta}$, $\delta$ may
have $V$-cofinality larger than $\omega_1$. However we have overcome the problem since
we have defined the sets $K(\alpha,\delta,\beta)$ more carefully than in the previous section. 
\hfill$\blacksquare$\medskip

In fact the theorem can be proved under the milder  assumptions
that $V$ and $W$ have the same cardinals, the same reals, $W\models\MRP$ and, for
every cardinal $\kappa$, there is in $V$ a partition
$\{A_s:s\in\kappa^{<\omega}\}$ of the points of $\kappa^+$ of
countable $V$-cofinality into disjoint sets which are stationary in $W$. By a
recent result of Larson building on ideas of Todor\v{cevi\'c} it is known
that such partitions can be found in $\ZFC$  for $\kappa=\omega$
just assuming that $\omega_1^V=\omega_1^W$. It is open whether for
higher cardinals such partitions exists in $\ZFC$. A positive
answer to this question would entail that if $V\subseteq W$ are
models with the same reals and cardinals and $W\models \MRP$ then
$ORD^\omega \cap V=ORD^\omega \cap W$.

We also remark that what we need to run the proof of theorem \ref{mainth} is the following weak form of reflection for stationary sets: 
{\it
\begin{quote}
Let $\kappa$ be of countable cofinality. Then for every countable family $\{A_n:n\in\omega\}$ of stationary sets of points of countable cofinality of $\kappa^+$ there are $n$, $\delta$, $\beta$ and there are infinitely many $j$ such that $A_j\cap K(n,\delta,\beta)$ is not empty.
\end{quote}
}

It is possible to see that this {\it "Reflection Principle"} holds above a strongly compact cardinal. 
What we have shown in this paper is that this property holds also under $\MRP$. 
\bigskip

\noindent {\small {\bf Acknowledgements:} The results presented in
this paper are a part of the Ph.D. thesis which I am completing
under the direction of professors Alessandro Andretta of the
University of Torino and Boban Veli\v{c}kovi\'c of the University
Paris 7. Part of this research has been completed while visiting
the Erwin Schr\"odinger International Institute for Mathematical
Physics during the Fall of 2004 in the framework of the ESI Junior
Research Fellowship 2004. I wish to thank Boban Veli\v{c}kovi\'c
who suggested me to investigate the effects of $\PFA$ on the
arithmetic of singular cardinals and who is guiding my steps in
the study of forcing axioms. I thank also Andr\'es Caicedo and
Justin Moore
for many useful comments and suggestions.

\end{document}